\documentclass[12pt]{amsart}
\usepackage{pstricks}
 \textwidth=\textwidth \textheight=\textheight
\theoremstyle{plain}
\newtheorem{thm}{Theorem}

\newtheorem{lem}[thm]{Lemma}
\newtheorem{prop}[thm]{Proposition}
\newtheorem{corollary}[thm]{Corollary}

\theoremstyle{definition}

\date{August 26, 2003}

\newcommand{\set}[1]{ \left\{ #1 \right\} }
\newcommand{\ov}[1]{\overline{#1}}
\def\a{\alpha}
\def\b{\beta}
\def\st{\overline{st}}
\def\ww{{\mathbf w}}
\def\c{\mathbb C}
\def\ll{\lambda}
\def\tll{\tilde{\lambda}}
\def\d{\delta}
\def\dd{\gamma}
\def\p{{\mathcal P}}
\def\vv{{\mathbf v}}

\psset{unit=1pt} \psset{arrowsize=4pt 1} \psset{linewidth=1pt}
\newrgbcolor{mygreen}{.2 .7 .2}
\newrgbcolor{myred}{.7 .3 .2}
\newgray{mygray}{.90}
\countdef\x=23 \countdef\y=24 \countdef\z=25 \countdef\t=26

\def\hdom(#1,#2)#3{
\x=#1 \y=#2 \multiply\x by 16 \multiply\y by 16 \z=\x \t=\y
\advance\z by 32 \advance\t by 16
\psline(\x,\y)(\x,\t)(\z,\t)(\z,\y)(\x,\y) \advance\x by 16
\advance\y by 8 \rput(\x,\y){{\bf #3}}}

\def\vdom(#1,#2)#3{
\x=#1 \y=#2 \multiply\x by 16 \multiply\y by 16 \z=\x \t=\y
\advance\z by 16 \advance\t by 32
\psline(\x,\y)(\x,\t)(\z,\t)(\z,\y)(\x,\y) \advance\x by 8
\advance\y by 16 \rput(\x,\y){{\bf #3}}}

\def\rec(#1,#2,#3,#4){
\psline(#1,#2)(#3,#2)(#3,#4)(#1,#4)(#1,#2)
}

\begin{document}
\title{Growth diagrams, Domino insertion and Sign-imbalance}
\author{Thomas Lam}
\address{Department of Mathematics,
         M.I.T., Cambridge, MA 02139}

\email{thomasl@math.mit.edu}

\begin{abstract}
We study some properties of domino insertion, focusing on aspects
related to Fomin's growth diagrams \cite{Fom1,Fom2}.  We give a
self-contained proof of the semistandard domino-Schensted
correspondence given by Shimozono and White \cite{SW}, bypassing
the connections with mixed insertion entirely.  The correspondence
is extended to the case of a nonempty 2-core and we give two dual
domino-Schensted correspondences.  We use our results to settle
Stanley's `$2^{n/2}$' conjecture on sign-imbalance \cite{Sta} and
to generalise the domino generating series of Kirillov, Lascoux,
Leclerc and Thibon \cite{KLLT} .
\end{abstract}

\maketitle

\section{Introduction}
Recently in \cite{SW} Shimozono and White described a semistandard
generalisation of domino insertion giving a bijection between
colored biwords and pairs of semistandard domino tableaux of the
same shape. They connected domino insertion with Haiman's mixed
and left-right insertion algorithms \cite{Hai} and via this
obtained the semistandard analogue.  They also made explicit a
color-to-spin property of domino insertion.  This property appears
to have been used earlier by Kirillov, Lascoux, Leclerc and Thibon
\cite{KLLT} for some special colored involutions.

Earlier, van Leeuwen \cite{vL} had described domino insertion in
terms of Fomin's growth diagrams.  He connected Barbasch and
Vogan's left-right insertion description \cite{BV} with
Garfinkle's traditional bumping description \cite{Gar}.  He also
defines insertion in the presence of a 2-core.

Our first aim in this paper is to give a self contained proof of
the semistandard domino-Schensted correspondence, using elementary
growth diagram calculations to prove all the main properties of
the bijection which we also extend to the nonempty 2-core case.
Thus our approach allows us to avoid mention of mixed insertion
completely.  We also describe two dual domino-Schensted
bijections.  These are bijections between multiplicity free
colored biwords and pairs of semistandard domino tableaux which
have conjugate shapes.  Finally, we perform a detailed analysis of
symmetric growth diagrams for domino insertion.

The study of growth diagrams leads us to a number of applications.
These include a number of enumerative results for domino tableaux,
an application to sign-imbalance, and a collection of product
expansions for generating series of domino functions.

The sign $sign(T)$ of a standard Young tableaux $T$ is the sign of
its reading word.  The sign imbalance of a shape $\ll$ is defined
as \[ \sum_{SYT \, T: sh(T) = \ll} sign(T). \]
 That sign-imbalance is related to domino tableaux has been
made explicit in work of White \cite{Whi} and Stanley \cite{Sta}.
In particular, White gives a formula for the sign of the Young
tableaux $T(D)$ associated to a domino tableaux $D$:
\[
sign(T) = (-1)^{ev(D)}
\]
where $ev(D)$ is the number of vertical dominoes in even columns
of $D$. Domino tableaux are in bijection with hyperoctahedral
involutions and we prove that in fact $ev(D)$ is equal to the
number of barred two-cycles of $\pi$, where $D = P_d(\pi)$ is the
insertion tableaux of $\pi$.  This allows us to prove Stanley's
conjecture on sign-imbalance, our Theorem \ref{thm:sta}, which is
a 4-parameter generalisation of the following elegant result:
\[
\sum_{SYT \, T: sh(T) \vdash m} sign(T) = 2^{\lfloor m/2 \rfloor}.
\]

Carr\'e and Leclerc \cite{CL} and Kirillov, Lascoux, Leclerc and
Thibon \cite{KLLT} have studied certain generating functions
$H_\ll(X;q)$ for domino tableaux which we loosely call
\emph{domino functions}.  More general domino functions
$G_\ll(X;q)$ were developed also in \cite{LLT}, where they were
connected with the Fock space representation of $U_q(\hat{
{\mathfrak sl_2}})$. These are defined as
\[
G_\ll(X;q) = \sum_D q^{spin(D)}x^D
\]
where the sum is over all semistandard domino tableaux of shape
$D$.  The $H_\ll$ are defined by $H_\ll(X;q)=G_{2\ll}(X;q)$.
Product expansions of the sums $\sum H_\ll(X;q)$ and $\sum H_{\ll
\vee \ll}(X;q)$ were given in \cite{KLLT}.

By studying colored involutions we give a product expansion for a
3-parameter generalisation of the sum $\sum_\ll G_\ll(X;q)$. When
the parameters of this sum is specialised, we obtain both of the
product expansions of \cite{KLLT}.

\medskip

We now briefly describe the organisation of this paper.  In
Section \ref{sec:prelims} we give some notation and definitions
for domino tableaux and colored words.  We also give a description
of domino insertion bumping in an informal manner, following
mostly \cite{SW}.  In Section \ref{sec:growth}, we introduce and
study growth diagrams.  This is followed by a proof of the
semistandard domino-Schensted correspondence and a description of
the dual domino-Schensted correspondences.  The section ends with
a study of symmetric growth diagrams and some enumerative results.
In Section \ref{sec:sign}, we apply the results of Section
\ref{sec:growth} to sign-imbalance.  In Section \ref{sec:domfunc},
we combine the results of Section \ref{sec:growth} with a study of
the standardisation of colored involutions.  These lead to a
number of product expansions for generating series of domino
functions. In Section \ref{sec:ribbon}, we give some final remarks
concerning possible generalisations to longer ribbons.

\bigskip
{\bf Acknowledgements}  I am indebted to my advisor, Richard
Stanley, for introducing the subjects of domino tableaux and
sign-imbalance to me and for suggesting his conjecture for study.
My work on generating series of domino functions was inspired by
the sum $\sum (-1)^{ev(\ll)}\tilde{G}_\ll(X;-1)$ suggested to me
by him.
\bigskip

\section{Preliminaries}
\label{sec:prelims}
\subsection{Domino Tableaux}
We will let $[n] = \set{1,2,\ldots,n}$ throughout.

Let $\lambda = (\ll_1 \geq \ll_2 \geq \ldots \geq \ll_{l(\ll)} >
0)$ be a partition of $n$. We will often not distinguish between a
partition $\ll$ and its diagram (often called $D(\ll)$) but the
meaning will always be clear from the context.  The partition $\ll
\cup \mu$ is obtained by taking the union of the parts of $\ll$
and $\mu$ (and reordering to form a partition).  We denote by
$\tilde{\lambda}$ and $(\ll^{(0)},\ll^{(1)})$ the $2$-core and
$2$-quotient of $\lambda$ respectively (see \cite{Mac}).  Every
2-core has the shape of a staircase $\d_r = (r,r-1,\ldots,0)$ for
some integer $r \geq 0$. As usual, when $\ll$ and $\mu$ are
partitions satisfying $\mu \subset \ll$ we will use $\ll/\mu$ to
denote the shape corresponding to the set-difference of the
diagrams of $\ll$ and $\mu$.

We denote the set of partitions by $\p$ and the set of partitions
with 2-core $\d_r$ by $\p_r$.  The set of all partitions $\ll$
satisfying the two conditions:
\begin{align*}
 \tll &= \d_r \\ |\ll| &= \d_r + 2n
\end{align*}
will be denoted $\p_r(n)$.  Note that $\p = \cup_{r,n} \p_r(n)$.

\medskip

A (standard) domino tableaux  (SDT) $D$ of shape $\ll$ consists of
a tiling of the shape $\ll / \tll$ by dominoes and a filling of
each domino with an integer in $[n]$ so that the numbers are
increasing when read along either the rows or columns.  Here, $n$
is $\frac{1}{2}|\ll / \tll|$.  A domino is any $2 \times 1 $ or $1
\times 2$ shape, or equivalently, two adjacent squares sharing a
common edge.  The \emph{value} of a domino is the number written
inside it.  We will write $dom_i$ to indicate the domino with the
value $i$ inside.  We will also write $sh(D) = \ll$.  An
alternative way of describing a standard domino tableaux of shape
$\ll$ is by a sequence of partitions $\set{\tll=\ll^0 \subset
\ll^1 \subset \ldots \subset \ll^n = \ll}$, where $sh(dom_i) =
\ll_i/\ll_{i-1}$.

A semistandard domino tableaux (SSDT) $D$ of shape $\ll$ consists
of a tiling of the shape $\ll / \tll$ by dominoes and a filling of
each domino with an integer, so that the numbers are
non-decreasing when read along the rows and increasing when read
along the columns.  The \emph{weight} of such a tableaux $D$ is
the composition $wt(D) = (\mu_1, \mu_2, \ldots)$ where there are
$\mu_i$ occurrences of $i$'s in $D$. Let $v(D)$ be the number of
vertical dominoes in a domino tableaux $D$.  The spin $sp(D)$, is
defined as $v(D)/2$. The standardisation of a semistandard domino
tableaux $D$ of weight $\mu$ is a standard domino tableaux
$D^{st}$ obtained from $D$ by replacing the dominoes containing
$1$'s with $1,2,\ldots,\mu_1$ from left to right, the dominoes
containing $2$'s by $\mu_1 + 1, \mu_1+2, \ldots, \mu_1+\mu_2$, and
so on.

\begin{figure}[ht]
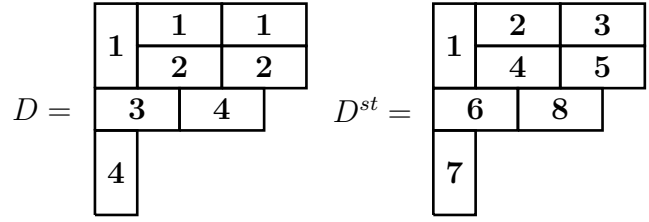

\pspicture(-20,0)(200,80) \rput(-20,40){$D=$} \vdom(0,0){4}
\hdom(0,2){3} \hdom(2,2){4} \vdom(0,3){1} \hdom(1,3){2}
\hdom(3,3){2} \hdom(1,4){1} \hdom(3,4){1}

\rput(105,40){$D^{st}=$} \vdom(8,0){7} \hdom(8,2){6}
\hdom(10,2){8} \vdom(8,3){1} \hdom(9,3){4} \hdom(11,3){5}
\hdom(9,4){2} \hdom(11,4){3}
\endpspicture
\caption{A domino tableaux $D$ with shape $(5,5,4,1,1)$ and weight
$(3,2,1,2)$ and its standardisation $D^{st}$.}
\label{fig:partition}
\end{figure}

More general skew (semi)standard domino tableaux are defined in a
similar manner.

\medskip

We should remark that Littlewood's 2-quotient map \cite{Lit} gives
a bijection between standard domino tableaux of shape $\ll$ and
pairs of standard Young tableaux of shapes $\ll^{(0)}$ and
$\ll^{(1)}$. This bijection generalises naturally to the
semistandard case. See for example \cite{CL}.  We will, however,
not be needing this bijection.

\subsection{Colored Words}
We will mostly follow the notation of \cite{SW} in this
subsection.

A \emph{letter} will be an integer with possibly a bar over it.

A \emph{colored word} is a word made of letters. A colored word
$\ww$ is a \emph{colored permutation} if each integer of $[n]$ is
used exactly once, for some $n$. Such a word will also be called a
hyperoctahedral permutation or a signed permutation. The set (in
fact group) of all such words will be denoted $B_n$. We define
$\ww^{neg}$ to be the word obtained from $\ww$ by converting all
the bars to negative signs. The weight of a word is defined in the
usual way, with the bars ignored.  The operation $ev$ removes the
bars from a colored word.  Thus if $\ww = (2\ov{3}\ov{1})$ then
$\ww^{ev} = (231)$.

A \emph{biletter} is an ordered pair of letters denoted ${x
\choose y}$.

A \emph{doubly colored biword} is a sequence of biletters ${x
\choose y}$ ordered canonically in the following way.  A biletter
${x \choose y}$ occurs before ${k \choose l}$ if and only if one
of the following holds:
\begin{enumerate}
\item $x < k$
\item $x = k$, both are unbarred, and $y^{neg} < l^{neg}$
\item $x = k$, both are barred and $l^{neg} < y^{neg}$.
\end{enumerate}

A doubly colored biword is a colored biword if only the bottom row
has bars.  For doubly colored biwords $\ww$, its inverse
$\ww^{inv}$ is obtained by swapping the top and bottom letters of
each biletter and then reordering.  For colored biwords $\ww$ its
inverse $\ww^{inv_r}$ is obtained by first moving the bars to the
top row and then performing the inverse normally.  Note that both
$inv$ and $inv_r$ are involutions.

The total color of a word or a colored word $\ww$, denoted
$tc(\ww)$, is the number of barred letters in the word.

Define $ev$ to be the operation which removes the bars from a
letter, word, biletter or a biword.

Standardisation $st$ is defined as follows for a colored biword
$\ww$. First define $\ww^{\st}$ by replacing the top row with
$\set{1,2,3,\ldots}$ from left to right when the biword is written
in order.  Then, define $\ww^{st}$ by (see \cite{SW})
\[
\ww^{st} = \ww^{\st \, inv \, \st \, inv}.
\]

For example, let $\ww$ be the colored biword
\[
\ww = \left( \begin{array}{ccccc}  1 & 1 & 2 & 3 & 3 \\ \ov{2} & 3
& 4 & \ov{1} & \ov{1} \\ \end{array} \right).
\]
Then $\ww$ has top weight $(2,2,1)$ and bottom weight $(2,1,1,1)$.
Its inverse $\ww^{inv_r}$ is given by
\[
\ww^{inv_r} = \left( \begin{array}{ccccc}  1 & 1 & 2 & 3 & 4 \\
\ov{3} & \ov{3} & \ov{1} & 1 & 2 \\ \end{array} \right).
\]
Its standardisation $\ww^{st}$ is computed as follows
\begin{align*}
\ww^{\st} &= \left( \begin{array}{ccccc}  1 & 2 & 3 & 4 & 5 \\
\ov{2} & 3 & 4 & \ov{1} & \ov{1} \\ \end{array} \right). \\
\ww^{\st \, inv} &= \left(\begin{array}{ccccc}  \ov{1} & \ov{1} &
\ov{2} & 3 & 4
\\ 5 & 4 & 1 & 2 & 3 \\ \end{array} \right) \\
\ww^{\st \, inv \, \st} &= \left(
\begin{array}{ccccc}  \ov{1} & \ov{2} & \ov{3} & 4 & 5 \\ 5 & 4 & 1 & 2 & 3 \\ \end{array} \right) \\
\ww^{\st \, inv \, \st \, inv} &= \left( \begin{array}{ccccc} 1 &
2 & 3 & 4 & 5
\\ \ov{3} & 4 & 5 & \ov{2} & \ov{1} \\ \end{array} \right)
\end{align*}

\begin{lem}
\label{lem:words} Let $\ww$ be a colored biword.  Then
\begin{align*}
\ww^{st} &= \ww^{inv \, \st \, inv \, \st} \\ \ww^{st} &= \ww^{\st
\, inv_r \, \st \, inv_r} \\ \ww^{inv_r \, st} &= \ww^{st \,
inv_r}.
\end{align*}
\end{lem}
\begin{proof}
The first statement is \cite[Proposition 12]{SW}.  The second
statement is essentially mentioned in \cite[Proposition 39]{SW}.
All three statements can be checked directly, which may be done
along the lines of Lemma \ref{lem:fixchange} (but is easier).
\end{proof}

We will occassionally identify a colored word $\ww$ or a
hyperoctahedral permutation $\pi$ with the associated colored
biword obtained by filling the top row with $\set{1,2,\ldots,n}$
from left to right.  In the latter case, $\pi^{inv}$ will be
identified with the lower row of the inverse of the resulting
biword.  This the usual inverse in the group $B_n$.

\subsection{Domino insertion}
\label{sec:domins} The normal Robinson-Schensted algorithm gives a
bijection between permutations of $S_n$ and pairs $(P,Q)$ of
standard Young tableaux (SYT) of size $n$ and the same shape.  A
semistandard generalisation of this was given by Knuth.  This is a
bijection between certain matrices with non-negative integer
entries (or alternatively unbarred biwords) and pairs of
semistandard Young Tableaux of the same shape.  We refer the
reader to \cite{EC2} for further details. Henceforth, familiarity
with usual Robinson-Schensted insertion will be assumed.

\medskip

In this section we describe the corresponding bijection for domino
tableaux in a traditional insertion `bumping' procedure. We will
follow the description given by Shimozono and White \cite{SW} for
the rest of this section where more details may be found.  As the
whole theory will be developed completely from the growth diagram
point of view in Section \ref{sec:growth}, we will not be
completely formal. The reader is referred to \cite{Gar},
\cite{SW}, \cite{vL} for full details.
\medskip

Let $D$ be a domino tableaux with $sh(D) = \ll$, no values
repeated, and $i$ a value which does not occur in $D$.  We will
describe how to insert both a vertical and horizontal domino with
value $i$ into $D$. Let $A \subset D$ be the sub-domino tableaux
containing values less than $i$.  If $\ll$ has a 2-core $\ll =
\tilde{\ll}$, then we will always assume that $\tilde{\ll} \subset
sh(A)$. We set $B$ to be the domino tableaux containing $A$ and an
additional vertical domino in the first column or an additional
horizontal domino in the first row labelled $i$.  Let $C = D/D'$
be the skew domino tableaux containing values greater than $i$.
Now we recursively define a bumping procedure as follows.

Let $(B,C)$ be a pair of domino tableaux (with no values repeated)
overlapping in at most a domino which contains the largest value
of $B$. The combined shape of $B$ and $C$ must be a valid skew
shape and the values of $C$ larger than those of $B$. Let $i < j$
be the largest and smallest values of $B$ and $C$ respectively.
Denote the corresponding dominoes by $\dd_i$ and $\dd_j$.  We now
distinguish four cases:
\begin{enumerate}
\item If $\dd_i$ and $\dd_j$ do not touch, then we set $B' =
B \cup \dd_j$ and $C' = C - \dd_j$.
\item If $\dd_i$ and $\dd_j$ intersect in exactly one
square, then we add a domino containing $j$ to $B$ so that the
shape of $B$ contains both $\dd_i$ and $\dd_j$ together with the
unique additional box which is diagonally outwards (right and
down) from $\dd_i \cap \dd_j$.  We set $C' = C - \dd_j$.
\item If $\dd_i = \dd_j$ and both are horizontal, then we
`bump' the domino $\dd_j$ to the next row, by setting $B'$ to be
the union of $B$ with an additional (horizontal) domino with value
$j$ one row below that of $\dd_i$.  We set $C' = C - \dd_j$.
\item If $\dd_i = \dd_j$ and both are vertical, then we
`bump' the domino $\dd_j$ to the next column, by setting $B'$ to
be the union of $B$ with an additional (vertical) domino with
value $j$ one column to the right of $\dd_i$.  We set $C' = C -
\dd_j$.
\end{enumerate}
This procedure is repeated with $(B,C)$ replaced by $(B',C')$
until the (skew) domino tableaux $C$ becomes empty.

The resulting $B$ tableaux will be denoted by $D \leftarrow i$ for
the insertion of a horizontal domino and $D \leftarrow
\overline{i}$ for a vertical domino.

\medskip

Let $\ww = w_1 w_2 \cdots w_n$ be a colored permutation and $\d_r$
be a 2-core assumed to be fixed throughout. Then the insertion
tableaux $P_d^r(\ww)$ is defined as $((\ldots((\d_r \leftarrow
w_1) \leftarrow w_2) \cdots ) \leftarrow w_n)$. The sequence of
shapes obtained in the process defines another standard domino
tableaux called the recording tableaux $Q_d^r(\ww)$.

As an example, the domino tableaux $P_d^0(\ov{3}42\ov{1})$ is
constructed as follows:
\begin{figure}[ht]
\pspicture(0,20)(304,80)
 \vdom(0,3){3}

\vdom(4,3){3} \hdom(5,4){4}

\hdom(9,4){2} \hdom(9,3){3} \vdom (11,3){4}

\vdom(16,3){1} \vdom(17,3){2} \hdom (16,2){3} \vdom(18,3){4}

\endpspicture
\caption{Insertion of $\ww = \ov{3}42\ov{1}$ into $\emptyset$.}
\label{fig:insertion}
\end{figure}

The following theorem will be proven in Section \ref{sec:growth}.
\begin{thm}
\label{thm:bij} Fix $r \geq 0$.  The above algorithm defines a
bijection between signed permutations $\pi \in B_n$ and pairs of
domino tableaux $(P,Q)$ of the same shape $\ll \in \p_r(n)$. This
bijection satisfies the equality
\[
tc(\pi) = sp(P_d(\pi)) + sp(Q_d(\pi)).
\]
\end{thm}

The insertion algorithm is due to Barbasch and Vogan \cite{BV} in
a different form (left-right insertion and jeu-de-taquin).  The
insertion described here in terms of bumping is essentially that
of Garfinkle \cite{Gar}. Van Leeuwen \cite{vL} proves that the
Barbasch-Vogan algorithm is the same as the bumping description,
and also shows that the bijection holds in the presence of a
2-core. That this algorithm sends total color to the sum of spins
seems to have been first used by Kirillov, Lascoux, Leclerc and
Thibon in \cite{KLLT} for certain hyperoctahedral involutions,
though no details or proofs are present.  More recently, the
color-to-spin property is made explicit by Shimozono and White in
\cite{SW}.

\par

Shimozono and White \cite{SW} only prove the color-to-spin
property in the absence of a 2-core.  However, the color-to-spin
property is proven by studying the spin change for all the `bumps'
in the insertion and these are unaffected by the presence of a
2-core. Thus the generalisation of the domino insertion bijection
to the 2-core case is immediate.  Shimozono and White also give a
semistandard generalisation of this bijection which is the case $r
= 0$ of the following theorem.  Their theory of domino insertion
is developed in conjunction with other combinatorial algorithms
including Haiman's mixed insertion and left-right insertion.

\begin{thm}
\label{thm:sembij} Fix a 2-core $\d_r$. There is a bijection
between colored biwords $\ww$ of length $n$ and pairs
$(P_d^r(\ww),Q_d^r(\ww))$ of semistandard domino tableaux with the
same shape $\ll \in \p_r(n)$ with the following properties:
\begin{enumerate}
\item
The bijection has the color-to-spin property:
\[
tc(\ww) = sp(P_d^r(\ww)) + sp(Q_d^r(\ww)).
\]
\item
The weight of $P_d^r(\ww)$ is the weight of the lower word of
$\ww$. The weight of $Q_d^r(\ww)$ is the weight of the lower word
of $\ww$.
\item
The bijection commutes with standardisation in the following
sense:
\begin{align*}
P_d^r(\ww)^{st} &= P_d^r(\ww^{st}). \\ Q_d^r(\ww)^{st} &=
Q_d^r(\ww^{st}).
\end{align*}
\end{enumerate}
\end{thm}
The proof of this will be left until the next section, where we
give an alternative description of domino insertion in terms of
growth diagrams.

\section{Growth Diagrams and Domino Insertion}
\label{sec:growth}
\subsection{Properties of Growth Diagrams}
The insertion algorithm of subsection \ref{sec:domins} can also be
phrased in terms of Fomin's growth diagrams \cite{Fom1,Fom2} (also
known as the poset-theoretic description, or language of shapes).
This was first made explicit by van Leeuwen \cite{vL}.  We will
show how growth diagrams are relevant to the semistandard
generalisation of domino insertion of \cite{SW}. Thus our aim will
be to give a short, stand-alone proof of Theorem \ref{thm:sembij}
using elementary considerations of growth diagrams only, bypassing
the connection with mixed insertion used by Shimozono and White.
Thus their lemma \cite[Lemma 33]{SW} is replaced by our Lemma
\ref{lem:ascent}. The use of growth diagrams make the
generalisation to the case of nonempty 2-core immediate.  In fact
one could use growth diagrams to define the entire correspondence
and develop the theory beginning from that.

\medskip

Let $M(i,j)$ be a $n \times n$ matrix taking values from
$\set{0,1,-1}$ thought of as the matrix representing a
hyperoctahedral permutation.  Thus it has one non-zero value in
each row or column.  We will take the row and column indices to
lie in $[n]$.

The growth diagram (of $M(i,j)$) is an array of partitions
$\ll_{(i,j)}$ for $1 \leq i,j \leq n+1$.  Two `adjacent'
partitions $\ll_{(i,j)}$ and $\ll_{(i+1,j)}$ or $\ll_{(i,j)}$ and
$\ll_{(i,j+1)}$ are either identical or differ by exactly one
domino.  Initially, all the partitions $\ll_{(1,j)}$ and
$\ll_{(i,1)}$ are set to the same partition $\mu$. For our
purposes this will usually be a partition satisfying $\mu =
\tilde{\mu}$. The remainder of the growth diagram will be
determined from $\mu$ and the data $M(i,j)$ according to the
following local rules.

Let $\ll = \ll_{(i,j)}$, $\mu = \ll_{(i+1,j)}$, $\nu =
\ll_{(i,j+1)}$, $\rho = \ll_{(i+1,j+1)}$ be the corners of a
`square'. Assume (inductively) that $\ll, \mu$ and $\nu$ are
known.  Then $\rho$ is determined as follows:
\begin{enumerate}
\item
If $M(i,j) = 1$ then it must be the case that $\ll = \mu = \nu$.
Obtain $\rho$ from $\ll$ by adding two to the first row.
\item
If $M(i,j) = -1$ then it must be the case that $\ll = \mu = \nu$.
Obtain $\rho$ from $\ll$ by adding two to the first column.
\item
If $M(i,j) = 0$ and $\ll = \mu$ or $\ll = \nu$ (or both) then
$\rho$ is set to the largest of the three partitions.
\item
Otherwise $M(i,j) = 0$ and $\nu$ and $\mu$ differ from $\ll$ by
dominoes $\dd$ and $\dd'$.  If $\dd$ and $\dd'$ do not intersect
then $\rho$ is set to be the union $\ll \cup \dd \cup \dd'$. If
$\dd \cap \dd'$ is a single square $(k,l)$, then $\rho$ is the
union of $\ll \cup \dd \cup \dd' \cup (k+1,l+1)$.  If $\dd = \dd'$
is a vertical domino then $\rho$ is obtained from $\ll \cup \dd$
by adding two to the column immediately to the right of $\dd$. If
$\dd = \dd'$ is a horizontal domino then $\rho$ is obtained from
$\ll \cup \dd$ by adding two to the row immediately below $\dd$.
\end{enumerate}
We will call these rules the \emph{local rules} of the growth
diagram.

\begin{prop}
\label{prop:growth} The above algorithm is well defined.  The
growth diagram models the insertion of the colored permutation
$\pi$ corresponding to $M(i,j)$ into a 2-core $\d_r$ (in fact more
generally any initial partition) .

The partition $\ll_{(i,j)}$ is the shape of the tableaux obtained
after the first $i$ insertions and restricted to values less than
$j$.  Thus $\set{\ll_{(n+1,j)}: j \in [n+1]}$ is a chain of
partitions determining $P_d^r(\pi)$ and $\set{\ll_{(i,n+1)}:i \in
[n+1]}$ is a chain of partitions determining $Q_d^r(\pi)$.
\end{prop}
\begin{proof}
This is proven via induction, by comparing domino insertion with
the local rules of the growth diagram.  The details can be found
in \cite{vL}.
\end{proof}

For example, Figure \ref{fig:growth} is the growth diagram
corresponding to the insertion procedure of Figure
\ref{fig:insertion}.

\begin{figure}[ht]
\pspicture(0,0)(320,320)

\psline[linecolor=gray, linewidth=0.5pt](0,0)(280,0)
\psline[linecolor=gray, linewidth=0.5pt](0,70)(280,70)
\psline[linecolor=gray, linewidth=0.5pt](0,140)(280,140)
\psline[linecolor=gray, linewidth=0.5pt](0,210)(280,210)
\psline[linecolor=gray, linewidth=0.5pt](0,280)(280,280)

\psline[linecolor=gray, linewidth=0.5pt](0,0)(0,280)
\psline[linecolor=gray, linewidth=0.5pt](70,0)(70,280)
\psline[linecolor=gray, linewidth=0.5pt](140,0)(140,280)
\psline[linecolor=gray, linewidth=0.5pt](210,0)(210,280)
\psline[linecolor=gray, linewidth=0.5pt](280,0)(280,280)

\rput(10,10){$\emptyset$} \rput(10,80){$\emptyset$}
\rput(10,150){$\emptyset$} \rput(10,220){$\emptyset$}
\rput(10,290){$\emptyset$} \rput(80,10){$\emptyset$}
\rput(80,80){$\emptyset$} \rput(80,150){$\emptyset$}
\rput(150,10){$\emptyset$} \rput(150,80){$\emptyset$}
\rput(150,150){$\emptyset$} \rput(220,10){$\emptyset$}
\rput(220,80){$\emptyset$} \rput(290,10){$\emptyset$}

\psline(75,215)(75,235)(85,235)(85,215)(75,215)
\psline(75,285)(75,305)(85,305)(85,285)(75,285)
\psline(145,215)(145,235)(155,235)(155,215)(145,215)
\psline(285,75)(285,95)(295,95)(295,75)(285,75)

\rec(215,145,235,155) \rec(285,145,305,165) \rec(215,215,235,235)
\rec(215,285,245,305) \rec(285,215,305,245)

\psline(145,285)(155,285)(155,295)(175,295)(175,305)(145,305)(145,285)
\psline(285,285)(305,285)(305,295)(315,295)(315,315)(285,315)(285,285)

\endpspicture
\caption{Growth diagram for the insertion of $\ww =
\ov{3}42\ov{1}$ into $\emptyset$.} \label{fig:growth}
\end{figure}

\begin{lem}
The local rules of a growth diagram are reversible in the
following sense.  Let $\ll = \ll_{(i,j)}$, $\mu = \ll_{(i+1,j)}$,
$\nu = \ll_{(i,j+1)}$, $\rho = \ll_{(i+1,j+1)}$ be the corners of
a `square' of the growth diagram.  Then $\rho$, $\mu$ and $\nu$
determine $\ll$ and $M(i,j)$.
\end{lem}
\begin{proof}
This is a simple verification of the local rules.
\end{proof}

Note, that there can be two legitimate standard domino tableaux
corresponding to $\set{\ll_{(i,n+1)}: i \in [n+1]}$ and
$\set{\ll_{(n+1,j)}: j \in [n+1]}$ which do not give a growth
diagram corresponding to an insertion procedure. For example if
$\ll_{(1,2)} = (2) = \ll_{(2,1)}$ and $\ll_{(2,2)}=(2,2)$ then
$\ll_{(1,1)}$ must be $\emptyset$.  This is not a valid growth
diagram for insertion as $\ll_{(1,1)} \neq \ll_{(2,1)}$.

\begin{lem}
\label{lem:rev} The correspondence
\[
\pi \rightarrow (P_d^r(\pi),Q_d^r(\pi))
\]
is a bijection between $\pi \in B_n$ and pairs of standard domino
tableaux of the same shape $\ll \in \p_r(n)$.
\end{lem}
\begin{proof}
The previous Lemma implies that this correspondence is injective.
As no dominoes can be removed from $\d_r$, the `initial' row and
column of the growth diagram ($\ll_{(1,j)}$ and $\ll_{(i,1)}$)
will consist completely of partitions equal to $\d_r$.  Thus
setting $\ll_{(i,n+1)}: i \in [n+1]$ and $\ll_{(n+1,j)}: j \in
[n+1]$ to two tableaux of the shape $\ll \in \p_r(n)$ will give a
growth diagram corresponding to the insertion of some
hyperoctahedral permutation $\pi$.
\end{proof}

\begin{lem}
\label{lem:sym} Let $\pi$ be a hyperoctahedral permutation. Domino
insertion possesses the symmetry property
\[
P_d^r(\pi) = Q_d^r(\pi^{inv}).
\]
\end{lem}
\begin{proof}
This is a consequence of the fact that the growth diagram local
rules are symmetric.
\end{proof}

\begin{lem}
\label{lem:colortospin} Domino insertion for hyperoctahedral
permutations $\pi$ possesses the color-to-spin property:
\[
tc(\pi) = sp(P_d^r(\pi)) + sp(Q_d^r(\pi)).
\]
\end{lem}
\begin{proof}
Let $\ll = \ll_{(i,j)}$, $\mu = \ll_{(i+1,j)}$, $\nu =
\ll_{(i,j+1)}$, $\rho = \ll_{(i+1,j+1)}$ be the corners of a
square of the growth diagram.  Then the Lemma follows from the
observation that
\[
sp(\rho/\mu) + sp(\rho/\nu) = sp(\mu/\ll) + sp(\nu/\ll) +
\begin{cases} 1 & \mbox{if $M(i,j) = -1$}\\
                0 & \mbox{otherwise}.
\end{cases}
\]
This can be checked by considering the local rules case by case.
\end{proof}

\begin{lem}
\label{lem:ascent} Let $\pi = \pi_1\cdots\pi_n$ be a colored
permutation. Then $\pi_i^{neg} < \pi_{i+1}^{neg}$ if and only if
$dom_i$ lies to the left of $dom_{i+1}$ in $Q^r_d(\pi)$.
\end{lem}
\begin{proof}
The main idea is to analyze a $1 \times 2$ rectangle of the growth
diagram.  Let $\ll_0 = \ll_{(i,j)}$, $\ll_1 = \ll_{(i+1,j)}$,
$\ll_2 = \ll_{(i+2,j)}$, $\mu_0 = \ll_{(i,j+1)}$, $\mu_1 =
\ll_{(i+1,j+1)}$ and $\mu_2 = \ll_{(i+2,j+1)}$ be the corners of a
$1 \times 2$ rectangle of the growth diagram.  We will call the
two squares of the $1\times 2$ rectangle the first and second
squares.  We further assume that $M(i,j) = M(i+1,j) = 0$.

Now suppose that $\a_0 = \ll_1/\ll_0$ and $\a_1 = \ll_2/\ll_1$ are
both dominoes so that $\a_0$ lies to the left of $\a_1$.  Then it
is easy to check that $\b_0 = \mu_1/\mu_0$ and $\b_1 =
\mu_2/\mu_1$ are both dominoes since $M(i,j) = M(i+1,j) = 0$.  We
claim that in fact $\b_0$ lies to the left of $\b_1$.  If $\ll_0 =
\mu_0$ this is trivial and most of the cases of the local rules
are a simple verification.

The only interesting case is when $\ll_1 = \mu_1$ and $\a_0$ is a
vertical domino.  In this case, $\b_0$ has moved to the right when
compared to $\a_0$. The key observation is that $\b_0$ is placed
in the column immediately to the right of $\a_0$, so it is either
still to the left of $\a_1$ or it overlaps $\a_1$.  When overlap
occurs, $\b_1$ will be moved further to the right and $\b_0$ will
remain to the left of $\b_1$.  This proves our claim.

\medskip

To show (one direction of) our lemma, we just need to check, case
by case, that the initial condition ($\a_0$ lying to the left of
$\a_1$) holds for $j = \max(\pi_i^{ev},\pi_{i+1}^{ev}) + 1$.  As
adding a new domino to the first column will be furthest to the
left, and adding a new domino to the first row will be the
furthest right this is a simple verification.  The claim implies
inductively that the same will continue to hold when we get to
$\ll_{(i,n+1)}$, $\ll_{(i+1,n+1)}$ and $\ll_{(i+2,n+1)}$, which
give exactly $dom_i$ and $dom_{i+1}$ of $Q_d^r(\pi)$.

The other direction of the lemma is proven in exactly the same
way, or one could replace `left' by 'above' and `row' by `column'.
\end{proof}

\begin{lem}
\label{lem:ascent2} Let $\pi = \pi_1\cdots\pi_n$ be a colored
permutation. Then $\left({\pi^{inv_r}}\right)_i^{neg} <
\left(\pi^{inv_r}\right)_{i+1}^{neg}$ if and only if $dom_i$ lies
to the left of $dom_{i+1}$ in $P^r_d(\pi)$.
\end{lem}
\begin{proof}
This is a consequence of Lemma \ref{lem:ascent} and Lemma
\ref{lem:sym}. \end{proof}

\medskip

We are now ready to prove the semistandard domino-Schensted
correspondence.
\begin{proof}[Proof of Theorem \ref{thm:sembij}]
That the correspondence exists for standard biwords is Lemma
\ref{lem:rev}.  Then the color-to-spin property follows from Lemma
\ref{lem:colortospin}.

\medskip

For the semistandard case, fix two weights $\mu$ and $\ll$ and let
these be the weights of the upper and lower words of a colored
biword $\ww$.  We define $P_d^r(\ww)$ by standardising the top row
first:
\[
P_d^r(\ww) = Q_d^r(\ww^{\st \, inv_r}).
\]
This is well defined since $\ww^{\st \, inv_r}$ has a colored
permutation as its lower word.  It is a semistandard domino
tableaux because of Lemma \ref{lem:ascent}.  This allows us to
define $Q_d^r(\ww)$ by
\[
Q_d^r(\ww) = P_d^r(\ww^{inv_r}).
\]

Next we show that these definitions commute with standardisation.
For example,
\begin{align*}
P_d^r(\ww)^{st} &= Q_d^r(\ww^{\st \, inv_r})^{st} \\
                &= Q_d^r(\ww^{\st \, inv_r \, \st}) \\
                &= Q_d^r(\ww^{st \, inv_r}) \\
                &= P_d^r(\ww^{st}).
\end{align*}
We have used Lemmas \ref{lem:words} and \ref{lem:sym}.  A similar
calculation proves that $Q_d^r(\ww)^{st} = Q_d^r(\ww^{st})$.

Since standardisation is injective (for both words and tableaux)
when the weights $\mu$ and $\ll$ are fixed, this proves that the
correspondence
\[
\ww \rightarrow (P_d^r(\ww),Q_d^r(\ww))
\]
is injective for colored biwords with fixed weights for the top
and bottom rows.  The color-to-spin property is also a consequence
of the standardisation procedure, as $tc(\ww) = tc(\ww^{st}) =
sp(P_d^r(\ww^{st}))+sp(Q_d^r(\ww^{st})) = sp(P_d^r(\ww)) +
sp(Q_d^r(\ww^{st}))$.

\medskip

Finally, one can show that correspondence is a surjection as
follows.  Suppose we are given a pair $(P,Q)$ of semistandard
domino tableaux of shape $sh(P)=sh(Q) \in \p_r(n)$ such that
$wt(P) = \ll$ and $wt(Q) = \mu$. Then we may obtain a colored
biword $\vv$ with standardised lower word by performing the
inverse correspondence (in the standard case) to
$(Q^{st},P^{st})$.  That the upper word can be converted to have
weight $\ll$ is a consequence of the `only if' part of Lemma
\ref{lem:ascent}. Thus $\vv$ satisfies $P_d^r(\vv) = Q^{st}$ and
$Q_d^r(\vv) = P$. Now perform the inverse correspondence to
$(P^{st},Q^{st})$, using Lemma \ref{lem:ascent} to prove that we
can change the upperword of $\vv^{inv_r}$ into weight $\mu$.

This completes the proof. \end{proof}

An alternative way of proving the surjectiveness of the
correspondence is by enumerating both colored words and pairs of
tableaux of the same shape.  Littlewood's 2-quotient map will
accomplish the latter.

For the case $r = 0$, it is easy to see that the definition used
in the proof agrees with that of Shimozono and White \cite{SW}.

\begin{corollary}
\label{cor:sym} The semistandard domino correspondence possesses
the symmetry property:
\[
P_d^r(\ww) = Q_d^r(\ww^{inv_r}).
\]
\end{corollary}
\begin{proof}
This is a consequence of the definition used in the proof.
\end{proof}

\subsection{Dual domino-Schensted correspondence}
In this section we give a description of two closely related dual
domino-Schensted correspondences.  They are bijections between
certain words and pairs of tableaux of the same shape, one of
which is semistandard and the other is column-semistandard.  For a
description of the dual RSK correspondence for Young tableaux see
\cite{EC2}.

A domino tableaux $D$ is \emph{column-semistandard} if its
transpose is semistandard.

A \emph{dual colored biword} is a colored biword such that the top
row is ordered as usual, but when the bottom row is used to order
two biletters, the reverse ordering is chosen.  Thus ${x \choose
y}$ precedes ${k \choose l}$ if
\begin{enumerate}
\item $x < k$, or
\item $x = k$ and $y^{neg} > l^{neg}$.
\end{enumerate}

The operator $\st$ is defined for dual colored biwords as usual by
standardising the top row.  The operation $inv_d$ changes dual
colored biwords to colored biwords and vice versa.  It swaps the
two letters of each biletter, moving the bar to the lower letter
if needed, and orders the biletters accordingly.

A colored biword or dual colored biword is called
multiplicity-free if any biletter ${i \choose j}$ occurs at most
once.  The same numbers may appear up to twice, but one must be
barred and the other non-barred.  For multiplicity-free biwords we
define the following new standardisation operation $std$ by
\[
\ww^{std} = \ww^{\st \, inv_d \, \st \, inv_d}.
\]

\begin{lem}
Let $\ww$ be a multiplicity free dual colored biword or colored
biword.  Then
\begin{align*}
\ww^{std} &= \ww^{inv_d \, \st \, inv_d \, \st} \\ \ww^{std \,
inv_d} &= \ww^{inv_d \, std}.
\end{align*}
\end{lem}
\begin{proof}
The proof is a direct verification, and very similar to Lemma
\ref{lem:words}.
\end{proof}

We may now define the two dual domino-Schensted correspondences
$\a$ and $\b$.  Let $\ww$ be a multiplicity-free dual colored
biword. Then we define $Q_\a^r(\ww)$ via domino-Schensted applied
to $\ww^{\st}$ where $\ww^{\st}$ is now treated as a colored
biword. To see that $Q_\a^r(\ww)$ is a column-semistandard domino
tableaux, we use Lemma \ref{lem:ascent}.  Also define $P_\a^r(\ww)
= P_d^r(\vv)$, where $\vv$ is the lower word of $\ww$.

Now let $\ww$ be a multiplicity-free colored biword.  We define
the correspondence $\b$ in a similar way.  Set $Q_\b^r(\ww)$ to be
$Q_d^r(\ww^{inv_d \, \st \, inv_d})$. We define $P_\b^r(\ww)$ by
turning $P_d^r(\ww^{inv_d \, \st \, inv_d})$ into a
column-semistandard tableaux of the same weight as the lower word
of $\ww$.  That this is possible is a consequence of Lemma
\ref{lem:ascent2}.

Note that both correspondences agree with the usual domino
correspondence when applied to hyperoctahedral permutations.

\begin{thm}
\label{thm:dual} Let $r \geq 0$ be fixed.  The map $\a$
\[
\a: \ww \rightarrow (P_\a^r(\ww), Q_\a^r(\ww))
\]
is a weight preserving bijection between multiplicity-free dual
colored biwords $\ww$ of length $n$ and pairs of tableaux $(P,Q)$
of the same shape $\ll \in \p_r(n)$ such that $P$ is semistandard
and $Q$ is column-semistandard.

The map $\b$
\[
\b: \ww \rightarrow (P_\b^r(\ww), Q_\b^r(\ww))
\]
is a weight preserving bijection between multiplicity-free dual
colored biwords $\ww$ of length $n$ and pairs of tableaux $(P,Q)$
of the same shape $\ll \in \p_r(n)$ such that $P$ is
column-semistandard and $Q$ is semistandard.

These maps satisfy the following properties:
\begin{enumerate}
\item They commute with standardisation.  Thus
\[
(P_\a^r(\ww)^{st}, Q_\a^r(\ww)^{st}) =
(P_d^r(\ww^{std}),Q_d^r(\ww^{std}))
\]
and similarly for $\b$.

\item The maps $\a$ and $\b$ are related by
\[
(Q_\a^r(\ww), P_\a^r(\ww)) =
(P_\b^r(\ww^{inv_d}),Q_\b^r(\ww^{inv_d})).
\]

\item Both maps have the color-to-spin property.
\end{enumerate}
\end{thm}
\begin{proof}
The proof is analogous to that of Theorem \ref{thm:sembij},
requiring use of Lemmas \ref{lem:ascent} and \ref{lem:ascent2}.
\end{proof}

\subsection{Statistics on Domino Tableaux}
In this subsection we will introduce and study a number of
statistics on partitions and domino tableaux.  Let $\ll$ be a
partition with 2-core $\tll$. Let $o(\ll)$ be the number of odd
rows of $\ll$ . Thus $o(\ll')$ is the number of odd columns. Let
\[
d(\ll) = \sum_{i=1}^{l(\ll/2)} \left\lfloor \frac{\ll_{2i}}{2}
\right\rfloor.
\]
Note that $d(\ll) = d(\ll')$ (see for example \cite{Sta}). Also
let
\[
v(\ll) = \sum_{i=1}^{l(\ll)} \left\lfloor \frac{\ll_i}{2}
\right\rfloor.
\]

Now let $D$ be a domino tableaux of shape $\ll$.  As before $v(D)$
is the number of vertical dominoes in $D$ and $sp(D) = v(D)/2$.
Let $ov(D)$ and $ev(D)$ be the number of vertical dominoes in odd
and even columns respectively.  Thus $sp(D) = (ov(D) + ev(D))/2$.
Let $mspin(\ll)$ be the maximum spin over all domino tableaux of
shape $\ll$. Similarly, let $ov(\ll)$ be the maximum of $ov(D)$
over all domino tableau of shape $\ll$.  Define $ev(\ll)$
similarly. The cospin of a domino tableaux $D$ is $cosp(D) =
mspin(\ll) - sp(D)$ (and is always an integer).

The following lemma is a strengthening of a lemma in \cite{Whi}.
\begin{lem}
\label{lem:ovev} Let $D$ be a domino tableaux of shape $\ll$ with
2-core $\tll$. Then
\begin{equation}
\label{eq:ovev} ov(D) - ev(D) = \frac{o(\ll) - o(\tll)}{2}.
\end{equation}
\end{lem}
\begin{proof}
We proceed by induction on the size of $\ll$, while keeping $\tll$
fixed. When $D$ has shape $\tll$ then both sides are 0.  Now let
$D$ have shape $\ll$ and suppose the Lemma is true for all shapes
$\mu$ that can be obtained from $\ll$ by removing a domino. Let
$\dd$ be the domino with the largest value in $D$.  Removing $\dd$
from $D$ gives a domino tableaux $D'$ for which (\ref{eq:ovev})
holds.  If $\dd$ is a horizontal domino then neither side changes.
If $\dd$ is a vertical domino in an odd row then both sides
decrease by 1 (changing from $D$ to $D'$).  If $\dd$ is a vertical
domino in an even row then both sides increase by 1.
\end{proof}

Note that this implies that a domino tableaux $D$ which has the
maximum spin (amongst all domino tableaux of shape $\ll$) will
also have the most number of odd vertical and even vertical
dominoes.  Thus for example, $mspin(\ll) = ev(\ll) + ov(\ll)$.

\subsection{Symmetric Growth Diagrams}
We now specialise to the case where the matrix $M_{\pi}(i,j)$
corresponds to a hyperoctahedral involution $\pi$.  Thus
$M_{\pi}(i,j)$ is symmetric and $\pi$ satisfies $\pi^2 = 1$ in the
group $B_n$. The hyperoctahedral involution $\pi$ will consist of
a number of fixed points, barred fixed points, two-cycles and
barred two-cycles.  For example, let $\pi =
(1\ov{6}\ov{3}54\ov{2}\ov{7})$. Then $\pi$ has one fixed point,
two barred fixed points, one two-cycle and one barred two-cycle.

\medskip

In this case we obtain the following proposition, part of which
was first observed by van Leeuwen \cite{vL}.

\begin{prop}
\label{prop:sym} Let $\pi \in B_n$ be a hyperoctahedral
involution. Suppose $\pi$ has $a$ fixed points, $b$ barred fixed
points, $c$ two-cycles and $d$ barred two-cycles.  Fix a 2-core
$\d_r$.  Let the insertion tableaux $P_d^r(\pi) = Q_d^r(\pi)$ of
$\pi$ into $\d_r$ have shape $\ll = sh(P_d^r(\pi))$ (which
satisfies $\tll = \d_r$). Then
\begin{align*}
sp(P_d^r(\pi)) &= \frac{b}{2} + d \\ \frac{o(\ll)-o(\d_r)}{2} &= b
\\ \frac{o(\ll') - o(\d_r)}{2} &=c \\ d(\ll)-d(\d_r) &= c + d.
\end{align*}
\end{prop}
\begin{proof}
Since $P_d^r(\pi) = Q_d^r(\pi)$ for a hyperoctahedral involution
by Lemma \ref{lem:sym}, the first equation is a consequence of the
color-to-spin property of Theorem \ref{thm:bij}.  For the other
statements, note that the symmetry of $M_{\pi}(i,j)$ and of the
local rules of the growth diagram imply that the growth diagram
$\ll_{(i,j)}$ itself is symmetric.  We focus our attention on the
partitions $\ll_{(i,i)}$. If $M_{\pi}(i,i) = 1$ then
$\ll_{(i+1,i+1)}$ has two boxes added to its first row, and so
$o(\ll_{(i+1,i+1)}') = o(\ll_{(i,i)}') + 2$. Similarly, if
$M_{\pi}(i,i) = -1$ then $o(\ll_{(i+1,i+1)}) = o(\ll_{(i,i)}) +
2$.  In both cases $d(\ll_{(i,i)}) = d(\ll_{(i+1,i+1)})$.
\par
If $M_{\pi}(i,i) = 0$ and $\ll_{(i+1,i)} = \ll_{(i,i)} =
\ll_{(i,i+1)}$ then $\ll_{(i,i)} = \ll_{(i+1,i+1)}$.  The only
remaining case is if $\ll_{(i+1,i)}$ differs from $\ll_{(i,i)}$ by
a domino, in which case $\ll_{(i,i+1)} = \ll_{(i+1,i)}$ as well.
This implies that $\ll_{(i+1,i+1)}$ differs from $\ll_{(i,i)}$ by
two dominoes in two adjacent columns or rows.  Regardless, the number
of odd columns and rows is unchanged while $d(\ll_{(i+1,i+1)}) =
d(\ll_{(i,i)}) + 1$.
\end{proof}

\begin{corollary}
\label{cor:ev} Let $D = P_d(\pi)$ correspond to a hyperoctahedral
involution $\pi$ with $b$ barred fixed points and $d$ barred
two-cycles. Then
\begin{align*}
ev(D) &= d. \\
ov(D) &= b+d.
\end{align*}
\end{corollary}
\begin{proof}
As before, let $\pi$ have $b$ barred fixed points.  Then by
Proposition \ref{prop:sym},
\[
ev(D) + ov(D) = 2sp(D) = b + 2d.
\]
Combining Lemma \ref{lem:ovev} with Proposition \ref{prop:sym}
again we have,
\[
ov(D) - ev(D) = \frac{o(\ll) - o(\tll)}{2} = b.
\]
Subtracting the two equations and dividing by two, we obtain the
first result.  Summing the two equations give the second result.
\end{proof}

The significance of this Corollary will become apparent in Section
\ref{sec:sign}.

\subsection{Some Enumeration for Domino Tableaux} Let $f^\ll$ be
the number of SYT of shape $\ll$.  The Robinson-Schensted
algorithm for standard Young tableaux (SYT) leads to a number of
enumerative results including the following well known result.

\begin{prop}
Let $n \geq 1$.  Then
\begin{equation}
\label{eq:nfact} \sum_{\ll \vdash n} (f^{\ll})^2 = n!.
\end{equation}
\begin{equation}
\label{eq:inv} \sum_{\ll \vdash n} f^{\ll} = t(n).
\end{equation}
\end{prop}

We can easily generalise these to domino tableaux.  Define
\[
d^{\ll}(q) = \sum_{SDT \, D: sh(D) = \ll} q^{spin(D)}.
\]
It is unlikely that a `hook-length' formula holds for
$d^{\ll}(q)$.  Note that $d^{\ll}(q)$ depends on more than just
the 2-quotient $(\ll^{(0)},\ll^{(1)})$ of $\ll$.  For example,
$(3,1,1)$ and $(2,2)$ have the same 2-quotient but $d^{(3,1,1)}(q)
= 2q^{1/2}$ and $d^{(2,2)}(q) = 1+q$.  A cospin version of
$d^{\ll}(q)$ for more general ribbon tableaux was studied by
Schilling, Shimozono and White in \cite{SSW}.

We have the following analogue of (\ref{eq:nfact}):
\begin{prop}
Let $n \geq 1$ and $r \geq 0$ be fixed.  Then
\[
\sum_{\ll}\left( d^{\ll}(q) \right )^2 = (1+q)^n n!
\]
where the sum is over all partitions $\ll \in \p_r(n)$.
\end{prop}
\begin{proof}
This is an immediate consequence of the bijection in Theorem
\ref{thm:bij}.
\end{proof}

Now define $h_r(n)$ as follows:
\[
h_r(n) = \sum_{\ll \in \p_r(n)}
a^{(o(\ll)-o(\d_r))/2}b^{(o(\ll')-o(\d_r))/2}c^{d(\ll)-d(\d_r)}d^{\ll}(q).
\]
When $a=b=c=q=1$, this is the number of hyperoctahedral
involutions in $B_n$ and thus a domino analogue of $t(n)$.

\begin{prop}
\label{prop:hypinv} The function $h(n) = h_r(n)$ does not depend
on $r$. It satisfies the recursion
\[
h(n+1) = (b + aq^{1/2})h(n) + nc(1+q)h(n-1).
\]
The exponential generating function defined as
\[
E_{h} = \sum h(n) \frac{t^n}{n!}
\]
is given by the formula
\[
E_{h} = exp \left((b+aq^{1/2})t + c(1+q)\frac{t^2}{2} \right).
\]
\end{prop}

\begin{proof}
That $h_r(n)$ does not depend on $r$ follows from the fact that
the tableaux being enumerated are in bijection with
hyperoctahedral involutions.  Furthermore, the bijection preserves
the appropriate weighting according to Proposition \ref{prop:sym}.
Thus we are in fact enumerating hyperoctahedral involutions.

The recursion for $h(n)$ is immediate from the construction of a
hyperoctahedral involution from barred and non-barred fixed points
and two-cycles.

For the exponential generating function, we can use the
exponential formula (see \cite[Corollary 5.1.6]{EC2}).  Thus we
think of a hyperoctahedral involution as a partition of $[n]$ into
one and two element subsets.  The one element subsets can be given
a weight of $b$ or $aq^{1/2}$ while the two element subsets can be
given a weight of $c$ or $cq$.
\end{proof}

\section{Sign-Imbalance and Stanley's Conjecture}
\label{sec:sign} Sign Imbalance can be defined for posets in
general, but we will only concern ourselves with the posets
arising from partitions.

Let $T$ be a standard Young tableaux.  Its reading word
$reading(T)$, for our purposes, will be obtained by reading the
first row from left to right, then the second row, and so on.  We
set $sign(T) = sign(reading(T))$ where $reading(T)$ is treated as
a permutation.

Let $\ll$ be a partition.  Then we set
\[
I_\ll = \sum_T sign(T)
\]
where the sum is over all standard Young tableaux $T$ of shape
$\ll$.  We say $I_\ll$ is the \emph{sign-imbalance} of $\ll$.

\medskip

It is not difficult to see that $I_\ll$ is related to domino
tableaux.  Suppose $\ll$ has no 2-core, then define an involution
on standard Young tableaux of shape $\ll$ by swapping $2i-1$ with
$2i$ for the smallest possible value of $i$ where this is
possible.  If no such swap is possible the tableaux is fixed by
the involution.

The fixed points correspond exactly to the standard domino
tableaux of shape $\ll$. We obtain a standard Young tableaux
$T(D)$ from a standard domino tableaux $D$, by filling the domino
with a 1 with the values 1 and 2, the domino with a 2, with the
values 3 and 4, and so on.

When $\ll$ has 2-core $\d_1$ (a single box) then we use an
involution which swaps $2i$ with $2i+1$ for the smallest value of
$i$ where it is possible.  Again, the fixed points are the
standard domino tableaux of shape $\ll$.

\par

It is easy to see that these involutions are sign-reversing on
tableaux which are not fixed points and thus we obtain the
following proposition.

\begin{prop}
Let $r \in \set{0,1}$, $n \geq 1$ and $\ll \in \p_r(n)$.  Then
\[
I_\ll = \sum_{sh(D) = \ll} sign(D)
\]
where the sum is over standard domino tableaux of shape $\ll$ and
the sign of a domino tableaux $D$ is the sign of the corresponding
standard Young tableaux $T(D)$ .
\end{prop}

For other values of $r$, we have the following result, see
\cite{Sta}:
\begin{prop}
Let $\ll$ have 2-core $\d_r$ for $r>1$, then \[I_\ll = 0.\]
\end{prop}

There is another natural involution on standard Young tableaux of
which standard domino tableaux are the fixed points.  This is
Sch\"{u}tzenberger's involution $S$, also known as evacuation. The
fixed points of this involution are exactly the domino tableaux of
shape $\ll$ satisfying $\tll = \d_r$ for $r \in \set{0,1}$ (see
\cite{vL}). For a fixed shape $\ll$, Stanley \cite{Sta} has shown
that $S$ is either always parity-reversing or parity-preserving.

\par
By analysing the positions of horizontal and vertical dominoes in
a standard domino tableaux, White \cite{Whi} proves the following
proposition.

\begin{prop}
\label{prop:white} Let $D$ be a domino tableaux of shape $\ll$
which has 2-core $\emptyset$ or $\d_1$.  Then
\[
sign(D) = (-1)^{ev(D)}.
\]
\end{prop}

White has also given an explicit formula (in terms of shifted
tableaux) for the sign-imbalance of partitions which have
`near-rectangular' shape.

Combining Proposition \ref{prop:white} with Corollary \ref{cor:ev}
we obtain the following theorem.
\begin{thm}
\label{thm:ev} Fix $r \in \set{0,1}$.  Let $\pi$ be a
hyperoctahedral involution. Then the sign of its insertion
tableaux $sign(P_d^r(\pi))$ is equal to the number of barred
2-cycles.
\end{thm}
\begin{proof}
This follows immediately from Corollary \ref{cor:ev} and
Proposition \ref{prop:white}.
\end{proof}

We can now prove the following conjecture of Stanley \cite{Sta},
known as the `$2^{\lfloor n/2 \rfloor}$' conjecture.

\begin{thm}
\label{thm:sta} Let $m \geq 1$ be an integer.  Then
\[
\sum_{\ll \vdash m}
x^{v(\ll)}y^{v(\ll')}q^{d(\ll)}t^{d(\ll')}I_\ll = (x+y)^{\lfloor
m/2 \rfloor}.
\]
Note that $d(\ll) = d(\ll')$ so that one of $q$ and $t$ is not
needed.
\end{thm}
\begin{proof}
Since $I_\ll = 0$ for $\ll$ with a 2-core larger than $\d_1$, we
may assume the sum is over $\ll \in \p_r(n)$, for the unique $r
\in \set{0,1}$ and $n$ satisfying $2n + r = m$.  Note that
$o(\d_1) = o(\d_1') = 1$ and $d(\d_1) = 0$.

The standard domino tableaux of such shape correspond exactly to
hyperoctahedral involutions $\pi \in B_n$.  We define an
involution $\alpha$ on all such $\pi$ by turning the two-cycle
$(i,j)$ with the smallest value of $i$ from barred to non-barred
or vice versa, if such an $i$ exists. By Theorem \ref{thm:ev},
$\alpha$ is sign-reversing for domino tableaux which are not fixed
points. Furthermore, by Proposition \ref{prop:sym}, all of the
statistics $o(\ll) - r$, $o(\ll') - r$ and $d(\ll)$ remain fixed
by $\alpha$.

The fixed points of $\alpha$ are exactly the hyperoctahedral
involutions without two-cycles.  Hence we obtain, using
Proposition \ref{prop:sym}
\[
\sum a^{(o(\ll)-r)/2}b^{(o(\ll')-r)/2}c^{d(\ll)}I_\ll = (a+b)^{n}.
\]

To change this into the form of Stanley's conjecture, observe that
$2v(\ll) + o(\ll) = m = 2n+r$ implying that $(o(\ll)-r)/2 = n -
v(\ll)$ and similarly for $v(\ll')$ and $o(\ll')$.  Now substitute
this and also $x=1/a$ and $y=1/b$.  Finally multiply both sides by
$(xy)^n$.
\end{proof}

Note that the fixed points of $\alpha$ in the proof are exactly
the domino tableaux which are hook shaped.  That these give the
right hand side of the conjecture was shown by Stanley \cite{Sta}.
When we set $x=y=q=1$ we obtain the following signed analogue of
(\ref{eq:inv}):
\[
\sum_{SYT \,T} sign(T) = 2^{\lfloor n/2 \rfloor}
\]
where the sum is over all standard Young tableaux $T$ of size $n$.

\section{Domino Generating Functions}
\label{sec:domfunc} Let $\Lambda$ denote the ring of symmetric
functions in a set of variables $X = (x_1,x_2,\ldots)$ taking
coefficients in $\c$ (though the coefficient field will not affect
the results). Its completion, $\tilde{\Lambda}$ includes symmetric
power series of unbounded degree (though the coefficient of a
monomial $m_\ll$ will always be well defined).

Carr\'{e} and Leclerc have defined symmetric functions
$H_\ll(X;q)$ via semistandard domino tableaux, in the same way
that Schur functions arise from semistandard Young tableaux.
Slightly more general functions $G_\ll(X;q)$ were used in
\cite{LLT} and the two are connected via $H_\ll(X;q) =
G_{2\ll}(X;q)$.

Let $\ll$ be a partition.  Define
\[
G_\ll = \sum_{D} q^{sp(D)}x^{wt(D)}
\]
where the sum is over all semistandard domino tableaux of shape
$\ll$ and $x^\mu := x_1^{\mu_1} x_2^{\mu_2} \ldots$ for a
partition $\mu$.  There is a cospin version of this function which
we will not need.  In the notation of \cite{LLT}, our $G_\ll$
would be denoted $G_{\ll/\tll}$.

That the $G_\ll$ are symmetric functions is a consequence of a
combinatorial interpretation of their expansion into Schur
functions given by Carr\'{e} and Leclerc.  As spin is not always
integral, the $G_\ll$ lie in the ring $\Lambda[q^{1/2}]$.  We will
call the $G_\ll$ \emph{domino functions}.  Theorem
\ref{thm:sembij} leads immediately to the following \emph{domino
Cauchy identity}.

\begin{prop}
Fix $r \geq 0$.  Then
\[
\sum_{\ll \in \p_r} G_\ll(X;q) G_\ll(Y;q) = \frac{1}{\prod_{i,j} (
1-x_iy_j)(1-qx_iy_j)}.
\]
\end{prop}

The dual domino-Schensted correspondence of Theorem \ref{thm:dual}
leads to the following \emph{dual domino Cauchy identity}.
\begin{prop}
Fix $r \geq 0$. Then
\[
\sum_{\ll \in \p_r} q^{|\ll/\d_r|/2}G_\ll(X;q) G_{\ll'}(Y;q^{-1})
= \prod_{i,j} ( 1+x_iy_j)(1+qx_iy_j).
\]
\end{prop}
\begin{proof}
This follows from the fact that a column-semistandard domino
tableaux $D$ with $m$ dominoes is a semistandard domino tableaux
$D'$ of the conjugate shape with spin given by
\[
sp(D') = \frac{m}{2} - sp(D).
\]
\end{proof}

\medskip

In \cite{KLLT}, Kirillov, Lascoux, Leclerc and Thibon give two
product expansions for certain sums of the $G_\ll$.  These will be
seen as specialisations of our Theorem \ref{thm:series}.  As the
paper \cite{KLLT} contains no proofs, our theorem can be
considered both as a proof and as a generalisation.

\medskip

We begin by studying closely the effect of standardisation on a
semistandard colored involution.

A colored word $\ww$ is said to be a colored involution if $\ww =
\ww^{inv_r}$.  Every such word is given by the number of fixed
points ${i \choose i}$, barred fixed points ${i \choose
\overline{i}}$, two-cycles ${i \choose j}...{j \choose i}$ and
barred two-cycles ${i \choose \overline{j}}...{j \choose
\overline{i}}$.  Let there be $a_i$, $b_i$, $c_{ij}$ and $d_{ij}$
of these respectively.  Thus $c_{ij} = c_{ji}$ and $d_{ij} =
d_{ji}$.

\begin{lem}
\label{lem:fixchange} Let $\ww$ be a colored involution.  Then its
standardisation $\ww^{st}$ is a signed involution with $a$ fixed
points, $b$ barred fixed points, $c$ two-cycles and $d$ barred
two-cycles, where:
\[
a = \sum_i a_i,
\]
\[
b = \sum_i b_i - 2\sum_i \left\lfloor\frac{b_i}{2}\right\rfloor
\]
\[
c = \sum_{i < j}c_{ij}
\]
\[
d = \sum_{i <j}d_{ij} +
\sum_i\left\lfloor\frac{b_i}{2}\right\rfloor.
\]
In other words, the only change that occurs is that of barred
fixed points becoming barred two-cycles.
\end{lem}
\begin{proof}
It is clear from Lemma \ref{lem:words} that $\ww^{st}$ is also an
involution.

Fix an integer $i$.  Then in the colored word $\mathbf{w}$, the
fixed points of the form ${i \choose i}$ have exactly
\[
A = \sum_{j < i} (a_j + b_j + c_{jk} + d_{jk}) + b_i + \sum_{k}
d_{ik} + \sum_{k < i} c_{ki}
\]
biletters in front.  If we look at $\ww^{\st \, inv}$ the same
formula holds using a different ordering for the top row. Thus
when we standardise and take inverse and standardise again, this
set of biletters will receive identical numbers for both the top
and bottom row, and will give us $a_i$ fixed points.

Now consider barred fixed points ${i \choose \overline{i}}$. There
are
\[
A = \sum_{j < i} (a_j + b_j + c_{jk} + d_{jk}) + \sum_{k>i} d_{ik}
\]
biletters in front.  Now let us consider what happens when we
standardise the top row and take the inverse.  We will similarly
get all (barred or otherwise) fixed points of $j < i$ in front and
so on. The only possible difference are the biletters involving
$i$. The fixed points clearly make no contribution. Since the
ordering for the lower letter is reversed when the upper letter is
barred, the biletters occuring in front are only those of the form
${\overline{i} \choose j}$ where $j > i$.  There are exactly
$\sum_{k>i} d_{ik}$ of these, thus the collection of barred fixed
points ${i \choose \overline{i}}$ will get the same set of numbers
for the upper and lower biletters.  However, individually, the
numbers assigned for the two rows will be reversals of each other
due to the ordering on the bottom row induced by the bars on the
upper row.

Now consider what happens to the collection of biletters of the
form ${i \choose j}$ and $i \neq j$.  We need only show that these
all become two-cycles when $\ww$ is standardised.  Since
$\ww^{st}$ is an involution we only need to check that these
biletters do not become fixed points.  Such a biletter has between
\[
A = \sum_{l < i} (a_l + b_l + c_{lk} + d_{lk}) + b_i + \sum_{k}
d_{ik} + \sum_{k < j} c_{ki}
\]
and
\[
B = \sum_{l < i} (a_l + b_l + c_{lk} + d_{lk}) + b_i + \sum_{k}
d_{ik} + \sum_{k < j} c_{ki} + c_{ij} - 1
\]
biletters in front.  After standardisation, exactly the same
formula holds with $i$ swapped with $j$.  We see that the top and
bottom letters will never get the same number via standardisation
(in fact if $i < j$ then $i$ will become a smaller number than
what $j$ becomes).

Exactly the same analysis holds for a biletter of the form ${i
\choose \overline{j}}$ and $i \neq j$.
\end{proof}

As an example, let $\ww$ be the colored involution
\[
\ww = \left( \begin{array}{ccccccccc} 1 & 1 & 2& 2& 2 &3 &3&4&5 \\
\ov{3} & 3& \ov{2}& \ov{2}& \ov{2}& \ov{1}& 1 &5&4 \\ \end{array}
\right)
\]
with 3 barred fixed points, 2 two-cycles and 1 barred two-cycle.
Then its standardisation
\[
\ww^{st} = \left( \begin{array}{ccccccccc} 1&2&3&4&5&6&7&8&9 \\
\ov{6}&7&\ov{5}&\ov{4}&\ov{3}&\ov{1}&2&9&8\\ \end{array} \right)
\]
has 1 barred fixed point, 2 two-cycle and 2 barred two-cycles.

\begin{thm}
\label{thm:series} Let $r \geq 0$ be fixed.  Let $S(X;a,b,c,q) \in
\tilde{\Lambda}[X][a,b,c,q^{1/2}]$ be the symmetric power series
\[
S(X;a,b,c,q^{1/2}) = \sum_{\ll \in \p_r}
a^{(o(\ll)-o(\d_r))/2}b^{(o(\ll')-o(\d_r))/2}c^{d(\ll)-d(\d_r)}G_\ll(X;q).
\]
Then $S(X;a,b,c,q^{1/2})$ does not depend on $r$ and has a product
formula given by
\[
\frac{\prod_i (1+aq^{1/2}x_i)}{\prod_i (1-bx_i) \prod_i
(1-cqx_i^2) \prod_{i<j} (1-cx_ix_j)\prod{i<j}(1-cqx_ix_j)}.
\]
\end{thm}

\begin{proof}
Semistandard domino tableaux are in one-to-one correspondence with
colored involutions by Theorem \ref{thm:sembij} and Corollary
\ref{cor:sym}.  If $\ww$ is a colored involution then the shape
and spin of $P_d^r(\ww)$ is that of $P_d^r(\ww^{st})$ and thus we
may use Proposition \ref{prop:sym} and Lemma \ref{lem:fixchange}
to calculate the contributions each colored involution makes to
the weights $o(\ll)$, $o(\ll')$, $d(\ll)$ and $sp(P_d^r(\ww))$.

Such colored involutions consist of a number of fixed points ${i
\choose i}$ corresponding to the product $\prod_i 1/(1-bx_i)$. The
barred fixed points ${i \choose \overline{i}}$ correspond to the
product $\prod_i (1+aq^{1/2}x_i)/(1-cqx_i^2)$ since according to
Lemma \ref{lem:fixchange} all but at most one of the barred fixed
points of each weight will pair to become a two-cycle upon
standardisation.  The two-cycles correspond to $\prod_{i<j}
1/(1-cx_ix_j)$ and the barred two-cycles correspond to
$\prod_{i<j} 1/(1-cqx_ix_j)$.
\end{proof}

There are a number of interesting specialisations.  We will set
$r=0$ for the next few examples.
\begin{enumerate}
\item
When $a=b=c=q^{1/2} = 1$, we obtain the square of a well known
identity:
\[
\left(\sum_{\ll \in \p} s_\ll(X)\right)^2 = \left(
\frac{1}{\prod_i (1-x_i) \prod_{i<j}(1-x_ix_j)} \right)^2.
\]
\item
Substituting $q^{1/2} = 0$ and using the fact that $G_{\ll}(X;0) =
s_\mu(X)$ for $\ll$ which satisfy $\ll = 2\mu$ (see \cite{CL}),
while $G_\ll(X;q) = 0$ for other $\ll \in \p_0$, we get
\[
\sum_{\ll \in \p} b^{o(\ll)}c^{v(\ll)}s_\ll(X) = \frac{1}{\prod_i
(1-bx_i) \prod_{i<j} (1-cx_ix_j)}.
\]
This is another well known identity which can be proved using
growth diagrams for normal RSK.
\item
The case $b=c=1$ and $a=0$ picks out the $G_\ll$ of the form
$G_{2\mu} = H_\mu$ and we obtain the first formula of \cite{KLLT}:
\[
\sum_{\ll} H_\ll(X;q) = \frac{1}{\prod_i (1-x_i)
\prod_{i<j}(1-x_ix_j) \prod_{i \leq j}(1-qx_ix_j)}.
\]
\item
The case $a=b=0$ and $c=1$ picks out the partitions of the form
$2\ll \vee 2\ll$ giving us the second formula of \cite{KLLT}:
\[
\sum_{\ll} H_{\ll \vee \ll}(X;q) =
\frac{1}{\prod_{i<j}(1-x_ix_j)\prod_{i \leq j}(1-qx_ix_j)}.
\]
\end{enumerate}

\medskip

Note that while $\sum G_\ll$ over $\ll \in \p_r(n)$ does not
depend on $r$,  the individual $G_\ll$ can differ greatly.  In
particular, two partitions $\ll$ and $\mu$ with the same
2-quotient but with $\tll \neq \tilde{\mu}$ may not have the same
$G$ function.  For example, $G_{(2,2)} = qs_2 + s_{1,1}$ while
$G_{(3,1,1)} = q^{1/2}(s_2+s_{1,1})$.  Both $(2,2)$ and $(3,1,1)$
have 2-quotient $\set{(1),(1)}$.

\section{Ribbon Tableaux}
\label{sec:ribbon} In this last section we make a few remarks
concerning which results might be generalised to ribbon tableaux.
We refer the reader to \cite{LLT} for the important definitions.

\medskip

Shimozono and White \cite{SW2} also give a spin-preserving
insertion algorithm for standard ribbon tableaux.  Unfortunately,
they stop short of giving a (spin-preserving) bijection between
words and pairs of semistandard tableaux.  Nevertheless, the
standard correspondence works.  It is a spin-preserving bijection
between pairs of standard ribbon tableaux and permutations $\pi$
of the wreath product $S_n\S C_p$. Again the involutions are in
bijection with standard ribbon tableaux and thus we obtain a
$p$-ribbon analogue of Proposition \ref{prop:hypinv} with an
identical proof.

\begin{prop}
Let $h(n)$ be the polynomial in $q$ defined as
\[
h(n) = \sum_{T} spin(T)
\]
where the sum is over all standard ribbon tableaux of size $n$
(and fixed $p$-core).  Then $h(n)$ satisfies the recurrence
\[
h(n+1) = (1+q^{1/2}+\ldots+q^{(p-1)/2})h(n) +
n(1+q+\ldots+q^{p-1})h(n-1)
\]
and has exponential generating function
\[
E_h(t) = exp\left( (1+q^{1/2}+\ldots+q^{(p-1)/2})t +
(1+q+\ldots+q^{p-1})\frac{t^2}{2} \right).
\]
\end{prop}

The statistics $o(\ll)$ and $d(\ll)$ are no longer suitable for
longer ribbons.  It seems likely that the statistic
\[
o_k(\ll) = \#\set{i: \ll_i \equiv k \mod p}
\]
may be interesting, but we have been unable to find any
applications.

\medskip
As Shimozono and White's ribbon correspondence can be phrased in
terms of growth diagrams, one might hope that a Lemma similar to
Lemma \ref{lem:ascent} can be shown in the same way -- this would
allow a semistandard ribbon correspondence to be developed.
Unfortunately this appears not to be the case, as ribbons may well
not `bump' to the next column or row but quite far away.  This
phenomenon occurs for certain longer ribbons regardless of whether
we insist upon column or row insertion/bumping.

\medskip

Possibly more promising is the following potential generalisation.
The sums over standard Young tableaux of size $n$
\begin{eqnarray*} \sum_T 1 &=& t(n) \\ \sum_T sign(T) &=&
2^{\lfloor n/2 \rfloor}
\end{eqnarray*} suggest that we might consider the sum
\[
\sum_T \chi(reading(T))
\]
for some other character $\chi$ of $S_n$.  If this were to be
related to $p$-ribbon tableaux and the wreath product $S_n \S C_p$
then $\chi$ should take $p^{th}$ roots of unity as its values. One
possibility is the (virtual) character which on the conjugacy
class of cycle type $\ll$ takes the value
\[
\chi(C_\ll) = \omega^{\ll - l(\ll)}
\]
for some $p^{th}$ root of unity $\omega$.

\end{document}